\documentclass{amsart}
\usepackage{amssymb}
\usepackage{amsthm}
\usepackage{amsmath}
\usepackage{amsxtra}
\usepackage{latexsym}
\usepackage{mathrsfs}
\usepackage[all]{xy}
\usepackage{graphicx}

\usepackage{url}

\newtheorem{thm}{Theorem}[section]
\newtheorem{lem}[thm]{Lemma}
\newtheorem{prop}[thm]{Proposition}
\newtheorem{cor}[thm]{Corollary}
\theoremstyle{definition}
\newtheorem{defn}[thm]{Definition}
\newtheorem{rmk}[thm]{Remark}
\newtheorem*{ack}{Acknowledgement}

\DeclareMathOperator{\id}{id}
\DeclareMathOperator{\Supp}{Supp}
\DeclareMathOperator{\Spec }{Spec}

\newcommand{\cquot}{\mathbin{\!/\!/\!}}

\title{On the Invariant Hilbert schemes
and Luna's {\'e}tale slice theorem}

\author{Yohsuke Matsuzawa}
\address{Graduate school of Mathematical Sciences, the University of Tokyo, Komaba, Tokyo, 153-8914, Japan}
\email{myohsuke@ms.u-tokyo.ac.jp}
\date{}

\begin{document}

\begin{abstract}
In this paper, we study local structures of invariant Hilbert schemes with Luna's {\'e}tale slice theorem.
We prove that in some cases the invariant Hilbert schemes are smooth
at a point which corresponds to a closed orbit.
\end{abstract}
\maketitle

\section{Introduction}

Let $k$ be an algebraically closed field of characteristic zero, 
$G$ a reductive algebraic group over $k$ and $X$  a scheme of finite type over $k$
with an action of $G$.
Let ${\rm Irr}(G)$ be the set of isomorphism classes of irreducible $G$-modules and 
$h\colon{\rm Irr}(G) \longrightarrow \mathbb{N}$ be a function (such a function is called a {\it Hilbert function}).
The {\it invariant Hilbert scheme} ${{\rm Hilb}^{G}_{h}}(X)$ is the scheme parametrizing $G$-stable
closed subschemes $Z$ of $X$ such that the coordinate ring of $Z$ is isomorphic to
${\bigoplus}_{M \in {\rm Irr}(G)}M^{\oplus h(M)}$ as a $G$-module.
The invariant Hilbert scheme has been introduced by Alexeev and Brion in \cite{ab05} for the case where $G$ is connected.
Its existence for non-connected $G$ is proved in \cite[Theorem 2.20]{brion}.

In \cite{ms10} Maclagan and Smith proved that the multigraded Hilbert scheme of the polynomial ring
$\mathbb{Z}[x,y]$ with a given Hilbert function is smooth and irreducible. Therefore the invariant Hilbert schemes are
smooth and irreducible when $G$ is diagonalizable and $X$ is the two-dimensional affine space
${\mathbb{A}}^{2}_{k}$ (in this case any action of it is isomorphic to a linear action).
When $G$ is a classical group and $X$ is a linear representation of $G$,
the smoothness of the invariant Hilbert scheme is studied by Terpereau and C. Lehn (\cite{ter12}, \cite{ter14a}, \cite{ter14b}
\cite{leter15}).
If the invariant Hilbert scheme is smooth, it gives a desingularization of the categorical
quotient $X\cquot G$.

In general, very little is known about the geometry of invariant Hilbert schemes.
However, several geometric properties such as smoothness and connectedness are expected 
when the problems are reduced to low dimensional cases.

In this paper, we study invariant Hilbert schemes by applying Luna's {\'e}tale slice theorem.
When $x\in X(k)$ is a $k$-valued smooth point and the orbit $Gx$ is closed, by using the slice theorem
we show that there are two {\'e}tale morphisms
\begin{align*}
& {\rm Hilb}^{G_{x}}_{h'}(V) \longrightarrow {\rm Hilb}^{G}_{h}(X),\\
 &{\rm Hilb}^{G_{x}}_{h'}(V) \longrightarrow {\rm Hilb}^{G_{x}}_{h'}(T_{x}V).
\end{align*}
Here $V$ is a slice at $x$, $G_x$ is the stabilizer of $x$,
$T_{x}V$ is the tangent space of $V$ at $x$, 
$h$ is a given Hilbert function with $h(0)=1$ (0 denotes the trivial representation)
and $h'$ is a suitable Hilbert function. In this paper, we study these morphisms
in detail (see Theorem \ref{b}, Remark \ref{red}).

If $Z$ is a $G$-stable closed subscheme of $X$ with underlying topological space
$Gx$ and $h$ is the Hilbert function of $Z$, then $h'$ has finite support (Theorem \ref{b}).
In this case ${{\rm Hilb}^{G_x}_{h'}}(T_xV)$ is an open subscheme of
${{\rm Hilb}_{n}}(T_xV)^{G_x}$ where $n=\sum_{M \in {\rm Irr}(G_x)}h'(M)$
and ${{\rm Hilb}_{n}}(T_xV)$ is the punctual Hilbert scheme.
Therefore, for example, ${{\rm Hilb}^{G}_{h}}(X)$ is smooth at $Z$ if
${{\rm Hilb}_{n}}(T_xV)$ is smooth.

In \S \ref{main thm} we state the main theorems, in \S \ref{proofs} we give the proofs
and in \S \ref{an example} we discuss an example of invariant Hilbert scheme of $ {\mathbb{A}}^{4}$ with an action of $ {\mathbb{G}}_{{\rm m}}^{2}$.

\section{The Main Theorems}\label{main thm}

Throughout this paper, $k$ is an algebraically closed field of characteristic zero and
$G$ is a reductive algebraic group over $k$.
Let ${\rm Irr}(G)$ be the set of isomorphism classes of irreducible representations of $G$.
The class of the trivial representation is denoted by $0$.
A variety over $k$ is a reduced separated scheme of finite type over $k$. 
A $G$-variety (resp. $G$-scheme) is a variety (resp. scheme) over $k$
with an action of $G$.
 
We use the following well known theorems (see \cite{luna} or \cite[Theorem 5.3, 5.4]{dr}).

\begin{prop}[Luna's {\'e}tale slice theorem]\label{slice}

Let $X$ be an affine $G$-variety. Let $x \in X(k)$ be a $k$-valued point such that 
the orbit $Gx$ is closed.
Then $G_x$, the stabilizer of $x$, is reductive by a theorem of Matsushima
and there exists a locally closed subvariety $V$ of $X$ such that     

\begin{enumerate}
\item[{\rm (1)}]  $V$ is affine and contains $x$;
\item[{\rm (2)}]  $V$ is $G_x$-stable;
\item[{\rm (3)}] the $G$-morphism $\psi \colon G {\times}^{G_x} V \longrightarrow X$ induced by 
the action of $G$ on $X$ is strongly {\'e}tale and the image $U$ of $\psi$ is a saturated open affine subset of $X$
(see Definition \ref{def} for the definition of strongly {\'e}tale morphisms and saturated open sets
and Definition \ref{def of induced scheme} for the definition of $G {\times}^{G_x} V$).
\end{enumerate}

We call such $V$ an {\'e}tale slice at $x$ for the action $G$ on $X$.
\end{prop}

When $X$ is smooth at $x$, we have the following.

\begin{prop}\label{smslice}
In the situation of Proposition \ref{slice}, if $X$ is smooth at $x$, $V$ can be chosen to be smooth and there exists 
a strongly {\'e}tale $G_x$-equivariant morphism $\varphi \colon V \longrightarrow T_{x}V$ such that
\begin{enumerate}
\item[{\rm (4)}] $\varphi (x)=0$, and $(d\varphi)_{x}=id$;
\item[{\rm (5)}]  $T_{x}X=T_{x}(Gx)\oplus T_{x}V$.
\end{enumerate}
\end{prop}

The {\'e}tale slice theorem can be used to analyze the local structure of quotients.
Since invariant Hilbert schemes can be regarded as "refined quotients",
the slice theorem is available to study the local structures
of invariant Hilbert schemes.

\begin{thm}\label{a}
 Let $h\colon{\rm Irr}(G) \longrightarrow \mathbb{N}$ be a Hilbert function with 
$h(0)=1$. In the situation of Proposition \ref{slice}, there exist finitely many Hilbert functions 
$h_{i} \colon {\rm Irr}(G_x) \longrightarrow \mathbb{N}, i=1,2,..,r$ with $h_{i}(0)=1$
and {\'e}tale morphisms of the invariant Hilbert schemes
  \[{{\rm Hilb}^{G_x}_{h_{i}}}(V) \longrightarrow {{\rm Hilb}^{G}_{h}}(X),\ i=1,2,..,r \]
 such that the union of  their images covers the open subset ${{\rm Hilb}^{G}_{h}}(U)$. 
 Here $U$ is the image of the map $G\times^{G_{x}}V \longrightarrow X$.

Moreover, if $X$ is smooth at $x$, then there are {\'e}tale morphisms

\[ {{\rm Hilb}^{G_x}_{h_{i}}}(V) \longrightarrow {{\rm Hilb}^{G_x}_{h_{i}}}({T_x}V),\ i=1,2,..,r .\]

\end{thm}

When $h(0)=1$, ${{\rm Hilb}_{h(0)}}(X\cquot G)$ is isomorphic to  the categorical quotient $X\cquot G$. 
Hence, via the Hilbert--Chow morphism
${{\rm Hilb}^{G}_{h}}(X) \longrightarrow {{\rm Hilb}_{h(0)}}(X\cquot G)$,
${{\rm Hilb}^{G}_{h}}(X)$ can be regarded as a refined quotient in this case.

At a point which corresponds to a closed orbit, the invariant Hilbert scheme is {\'e}tale locally
isomorphic to an invariant Hilbert scheme with finitely supported Hilbert function.

\begin{thm}\label{b}
Let $X$ be an affine $G$-variety. Let $x\in X(k)$ be a point such that 
$Gx$ is a closed orbit and $V$ an {\'e}tale slice at $x$. Let $Z\subset X$ be a $G$-stable closed 
subscheme such that the underlying topological space is $Gx$.
Let $h$ be the Hilbert function of $Z$. Suppose $h(0)=1$.
Then there is a finitely supported Hilbert function $h' \colon{\rm Irr}(G_x) \longrightarrow \mathbb{N}$ 
with $h'(0)=1$ and an {\'e}tale morphism
\[ \Psi \colon {{\rm Hilb}^{G_x}_{h'}}(V) \longrightarrow {{\rm Hilb}^{G}_{h}}(X) \]
such that the image of $\Psi$ contains $Z \in {{\rm Hilb}^{G}_{h}}(X)(k)$.

Moreover, if $X$ is smooth at $x$, then there is an {\'e}tale morphism
\[
{{\rm Hilb}^{G_x}_{h'}}(V) \longrightarrow {{\rm Hilb}^{G_x}_{h'}}({T_x}V).
\]
\end{thm}

In several cases, we can prove the invariant Hilbert scheme is smooth at a point
corresponding to a closed orbit by using Theorem \ref{b}.

\begin{cor}\label{c}
Let $X$ be an affine $G$-variety and $x\in X(k)$ a smooth point such that $Gx$ is a closed orbit.
Let $Z\subset X$ be a $G$-stable closed subscheme with Hilbert function $h$ such that the underlying topological space is $Gx$.
Suppose $h(0)=1$.
Let ${\rm dim}_xX $ denote the maximum of the dimensions of irreducible components of $X$ containing $x$. 
Then ${{\rm Hilb}^{G}_{h}}(X)$ is smooth at $Z$ if
\begin{enumerate}
\item[{\rm (1)}] ${\rm dim}_xX-{\rm dim}Z \leq 2$,
\item[{\rm (2)}] $Z$ is Gorenstein and ${\rm dim}_xX-{\rm dim}Z=3$, or
\item[{\rm (3)}] $Z$ is reduced.
\end{enumerate}

\end{cor}

As an obvious corollary, we have the following.

\begin{cor}
Let $X$ be a smooth affine $G$-variety such that all the orbits are closed.
{\rm(}In this case, the categorical quotient $X\cquot G$ is actually a geometric quotient that we denote $X/G$ for short.{\rm )}
Let $h\colon{\rm Irr}(G) \longrightarrow \mathbb{N}$ be a Hilbert function with $h(0)=1$.
\begin{enumerate}
\item[{\rm (1)}] If ${\rm dim}X/G \leq 2$, then ${{\rm Hilb}^{G}_{h}}(X)$ is smooth.
\item[{\rm (2)}] If ${\rm dim}X/G \leq 3$, then ${{\rm Hilb}^{G}_{h}}(X)$ is smooth at every point which represents a Gorenstein subscheme.
\item[{\rm (3)}] ${{\rm Hilb}^{G}_{h}}(X)$ is smooth at every point which represents a reduced subscheme.
\end{enumerate}
\end{cor}

\begin{proof}
Let $\pi \colon X \longrightarrow X/G$ be the quotient morphism.
Let $Z$ be a point of ${{\rm Hilb}^{G}_{h}}(X)(k)$.
Since $h(0)=1$, $\pi (Z)$ is one point. Therefore if we take a point $x\in Z(k)$,
then the underlying space of $Z$ is $Gx$.
Note that ${\rm dim}_xX-{\rm dim}Z \leq {\rm dim}X/G$.
Now the statements follow from Corollary \ref{c}.
\end{proof}

\section{Proofs}\label{proofs}

In order to prove Theorem \ref{a} and Theorem \ref{b}, we need some lemmas.

First, we investigate the relationship between the invariant Hilbert schemes of the source and target of the morphism in
Proposition \ref{slice} (3).
\begin{defn}\label{def}
Let $X,Y$ be affine $G$-schemes of finite type over $k$.
\begin{enumerate}
\item[{\rm (1)}] Let $\pi\colon X \longrightarrow X\cquot G$ be the categorical quotient.
An open subset $U \subset X$ is said to be saturated if
there is an open subset $V \subset X\cquot G$ such that ${\pi}^{-1}(V)=U$.
\item[{\rm (2)}]  A morphism $\varphi \colon X \longrightarrow Y$ is called strongly {\'e}tale if

\begin{enumerate}
\item[{\rm (i)}]  The following commutative diagram is cartesian
\[
\xymatrix{X \ar[r]^{\varphi} \ar[d] & Y \ar[d]
\\X\cquot G \ar[r]_{{\varphi}\cquot G} & Y\cquot G\ ;}
\]

\item[{\rm (ii)}]  the induced morphism $\varphi \cquot G$ is {\'e}tale.
\end{enumerate}
Note that $\varphi$ is {\'e}tale by the conditions (i) and (ii).
\end{enumerate}

\end{defn}

In Lemma \ref{satop} and \ref{stet}, we show that saturated open immersions and strongly {\'e}tale morphisms
induce morphisms between corresponding invariant Hilbert schemes.

\begin{lem}\label{satop}
Let $X$ be an affine $G$-scheme of finite type over $k$, $U \subset X$
be a saturated open affine subset. For every $h\colon{\rm Irr}(G) \longrightarrow \mathbb{N}$,
there is a canonical open immersion
\[ {{\rm Hilb}^{G}_{h}}(U) \longrightarrow {{\rm Hilb}^{G}_{h}}(X). \]
\end{lem}
\begin{proof}
We consider Hilbert schemes as contravariant functors from the category of noetherian $k$-schemes $({\rm Sch}/k)$
to the category of sets. We have to construct a natural transformation and prove it is
an open immersion.

Recall that a natural transformation $\eta \colon F \longrightarrow G$ of  contravariant functors
from $({\rm Sch}/k)$ to the category of sets is called an open immersion
when the following conditions are satisfied:
\begin{enumerate}
\item[\rm (1)]  For every $S \in ({\rm Sch}/k)$, $\eta(S) \colon F(S) \longrightarrow G(S)$ is injective.
\item[\rm (2)] For every $S \in ({\rm Sch}/k)$ and $g \in G(S)$, there exists an open subscheme $S_{g} \subset S$ such that
for any morphism $\varphi \colon T \longrightarrow S$ in $({\rm Sch}/k)$, the following are equivalent:
\begin{itemize}
\item $\varphi$ factors through $S_{g} \to S$.
\item $G(\varphi)(g)$ is contained in $\eta(T)(F(T))$.
\end{itemize}
\end{enumerate}

Now, let $S$ be a noetherian $k$-scheme and $W \subset U \times S$ be a $G$-stable closed subscheme
with Hilbert function $h$. Then we have the following commutative diagram.\\
\[
\xymatrix{&W \ar[d]^i\ar[dl]\\
W'\ar[r]_{i'}\ar[d] & U\times S \ar[r]^{j} \ar[d] & X\times S \ar[d]\\
W\cquot G \ar[r]_{i\cquot G} & U\cquot G\times S \ar[r]_{{j}\cquot G} & X\cquot G\times S}
\]
Here $i$ (resp. $j$) is the inclusion, $i\cquot G$ (resp. $j\cquot G$) is the induced closed immersion (resp. open immersion), $i'$ is
the base change of $i\cquot G$. Since $U$ is saturated, $j$ is the base change of $j\cquot G$.
Since $W$ is multiplicity finite (in the sense of \cite[Definition 2.13]{brion}), $W\cquot G$ is finite over $S$ and $j\cquot G \circ i\cquot G$
is finite. Therefore $j\circ i'$ is finite and so is $j\circ i$.
Thus $j\circ i$ is a finite immersion, therefore a closed immersion.
So we obtain an injective natural transformation ${{\rm Hilb}^{G}_{h}}(U) \longrightarrow {{\rm Hilb}^{G}_{h}}(X)$ by sending $i\colon W \longrightarrow U\times S$ to $j\circ i\colon W \longrightarrow X\times S$.

We check the condition (2).
Let $W\subset X\times S$ be a closed $G$-subscheme with Hilbert function $h$.
Let $Z=W \cap ((X\setminus U)\times S)$, where $X \setminus U$ is equipped with its reduced structure.
Then $Z$ is a multiplicity finite closed $G$-subscheme of $X \times S$.
Let $p\colon Z \longrightarrow S$ be the structure morphism. Then for any morphism $f\colon S'\longrightarrow S$,
\[
f^{*}W\in {{\rm Hilb}^{G}_{h}}(U)(S') \Longleftrightarrow f^{*}Z=\emptyset
\Longleftrightarrow f\ {\rm factors}\ S\setminus p(Z).
\]
But since $p\colon Z\longrightarrow Z\cquot G \longrightarrow S$ and $Z\cquot G$ is finite over $S$, $p(Z)$ is closed in $S$ and $S\setminus p(Z)$ is open;
hence ${{\rm Hilb}^{G}_{h}}(U)$ is an open subfunctor of ${{\rm Hilb}^{G}_{h}}(X)$. 
\end{proof}

\begin{lem}\label{stet}
Let $X, Y$ be affine $G$-schemes of finite type over $k$, $\varphi \colon  X \longrightarrow Y$
be a strongly {\'e}tale morphism. For any $h\colon {\rm Irr}(G) \longrightarrow \mathbb{N}$
with $h(0)=1$, $\varphi$ induces an {\'e}tale morphism
\[ \Phi \colon {{\rm Hilb}^{G}_{h}}(X) \longrightarrow {{\rm Hilb}^{G}_{h}}(Y). \]
If $\varphi$ is surjective, then $\Phi$ is surjective.
\end{lem}
\begin{proof}{\bf Step 1}.
First we construct a natural transformation from ${{\rm Hilb}^{G}_{h}}(X)$
to ${{\rm Hilb}^{G}_{h}}(Y)$.
Let $S$ be a noetherian $k$-scheme, and $W\subset X\times S$ be
a $G$-stable closed subscheme with Hilbert function $h$. Then we have the
following commutative diagram.
\[
\xymatrix{&W \ar[d]^i\ar[dl]_{j}\\
W'\ar[r]_{i'}\ar[d] & X\times S \ar[r]^{\varphi'=\varphi \times \id} \ar[d] & Y\times S \ar[d]\\
W\cquot G \ar[r]_{i\cquot G} & X\cquot G\times S \ar[r]_{{\varphi}\cquot G} & Y\cquot G\times S}
\]
Here $i$ stands for the inclusion, $i'$ is the base change of $i\cquot G$ and $j$ is the morphism
obtained from the universal property for the Cartesian square.
Since $\varphi$ is strongly {\'e}tale, $\varphi'$ is the base change of $\varphi \cquot G$.
By the assumption that $h(0)=1$, $W\cquot G \longrightarrow Y\cquot G\times S \longrightarrow S$ is an isomorphism
and $\varphi \cquot G \circ i\cquot G$ is a closed immersion.
Thus $\varphi' \circ i'$ is also a closed immersion since this is the base change of $\varphi \cquot G \circ i\cquot G$.
Hence $\varphi' \circ i=\varphi' \circ i' \circ j$ is a closed immersion. (Note that $j$ is a closed immersion since $i$ is.) 
Thus we get a natural transformation
\[
\Phi\colon {{\rm Hilb}^{G}_{h}}(X)(S) \longrightarrow {{\rm Hilb}^{G}_{h}}(Y)(S)\ ;(i\colon W\longrightarrow X\times S) \mapsto
(\varphi' \circ i \colon  W \longrightarrow Y\times S) .
\]

{\bf Step 2}.
Next we verify that the natural transformation $\Phi$ constructed in Step 1
is an {\'e}tale morphism. By the infinitesimal lifting property(c.f.\ \cite[Theorem 2.6.2]{fu}), we have to check that             
for an arbitrary noetherian affine scheme $\Spec A$ over ${{\rm Hilb}^{G}_{h}}(Y)$
and an ideal $I\subset A$ with $I^2 =0$,
\[
{\rm Hom}_{{{\rm Hilb}^{G}_{h}}(Y)}(\Spec A,{{\rm Hilb}^{G}_{h}}(X)) \longrightarrow
{\rm Hom}_{{{\rm Hilb}^{G}_{h}}(Y)}(\Spec A/I,{{\rm Hilb}^{G}_{h}}(X))
\]
is bijective.
To this end, we need to check that the map
\[
{\Phi}(A)^{-1}(Z) \longrightarrow {\Phi}(A/I)^{-1}(\tilde{Z}) 
\]
is bijective for every noetherian ring $A$ over $k$,
ideal $I\subset A$ with $I^2 =0$, and $Z \in {{\rm Hilb}^{G}_{h}}(Y)(A)$.
Here $\tilde{Z}$ is the image of $Z$ by
${{\rm Hilb}^{G}_{h}}(Y)(A) \longrightarrow {{\rm Hilb}^{G}_{h}}(Y)(A/I)$.
The situation is summarized in the following diagram:
\[
\xymatrix@C=36pt@R=5pt{
{{\rm Hilb}^{G}_{h}}(X)(A) \ar[r] \ar[ddd]_{\Phi (A)} & {{\rm Hilb}^{G}_{h}}(X)(A/I) \ar[ddd]^{\Phi (A/I)}\\
&\\
&\\
{{\rm Hilb}^{G}_{h}}(Y)(A) \ar[r] & {{\rm Hilb}^{G}_{h}}(Y)(A/I)\\
Z \ar@{}[u]|{\rotatebox[origin=c]{90}{$\in$}} \ar@{|->}[r] & \tilde{Z}. \ar@{}[u]|{\rotatebox[origin=c]{90}{$\in$}}
}
\]
Consider the following commutative diagram where every square is cartesian:
\[
\xymatrix{
\tilde{T} \ar[rd] \ar[ddd] \ar[rrr]^{\beta} & & &\tilde{Z} \ar[ld] \ar[ddd]\\
&{X\times \Spec A/I} \ar[r] \ar[d] &{Y\times \Spec A/I} \ar[d]&\\
& {X\times \Spec A} \ar[r] &{Y\times \Spec A} &\\
T \ar[ru] \ar[rrr]_{\alpha}& & &Z. \ar[lu]
}
\]
Then
\begin{align*}
{\Phi}(A)^{-1}(Z) \simeq \{ {\rm sections\ of}\ \alpha \} \simeq
\{ D\subset T\ {\rm open}\mid{\alpha}|_D\ {\rm is\ an\ isom}\};\\
{\Phi}(A/I)^{-1}(\tilde{Z}) \simeq \{ {\rm sections\ of}\ \beta \} \simeq
\{ \tilde{D}\subset \tilde{T}\ {\rm open}\mid{\beta}|_{\tilde{D}}\ {\rm is\ an\ isom}\}.
\end{align*}
Note that $Z$ is an affine scheme and the ideal $J$ of $\tilde{Z}$ in $Z$
satisfies $J^2=0$. Since $\alpha$ is {\'e}tale, these two sets are canonically bijective.

{\bf Step3}.
Suppose $\varphi$ is surjective. 
Let $Z\subset Y$ be a $G$-stable closed subscheme with Hilbert function $h$.
Then $Z\cquot G$ is one point and  it is denoted by $z\in Y\cquot G$.
Consider the following commutative diagram
\[
\xymatrix{
W \ar[rd] \ar[ddd] \ar[rrr] & & &Z \ar[ld] \ar[ddd]\\
&X \ar[r]^{\varphi} \ar[d] &Y \ar[d]&\\
& X\cquot G \ar[r]_{\varphi \cquot G} &Y\cquot G &\\
{(\varphi \cquot G)}^{-1}(z) \ar[ru] \ar[rrr]& & &\{z\}=\Spec k \ar[lu]
}
\]
where $W$ is the pull-back of $Z$ by $\varphi$. Since $\varphi \cquot G$ is {\'e}tale,
${(\varphi \cquot G)}^{-1}(z)$ is a disjoint union of $\Spec k$.
Therefore $W \longrightarrow Z$ has a section and $Z$ is in the image of
${{\rm Hilb}^{G}_{h}}(X)(k) \longrightarrow {{\rm Hilb}^{G}_{h}}(Y)(k)$.
\end{proof}

Next we investigate the relationship between the invariant Hilbert schemes of $G$ and $G_x$.

\begin{defn}\label{def of induced scheme}
Let $H \subset G$ be a reductive subgroup, 
and $X$ be an affine $H$-scheme of finite type over $k$.
Then $H$ acts on $G \times X$ by $h(g,x)=(gh^{-1},hx)$.
Then we denote the quotient $(G \times X)\cquot H$ by $G {\times}^H X$
and call it the induced scheme. 
We equip $G {\times}^H X$
with a $G$-scheme structure by $g(g',x)=(gg',x)$.
\end{defn}

\begin{lem}\label{ind}
Let $H \subset G$ be a reductive subgroup, 
and $X$ be an affine $H$-scheme of finite type.
Then for every $h\colon {\rm Irr}(G) \longrightarrow \mathbb{N}$,
there are finitely many Hilbert functions $h_{i} \colon {\rm Irr}(H) \longrightarrow \mathbb{N}$
with $h_{i}(0)=1$ such that
${{\rm Hilb}^{H}_{h_{i}}}(X)$ is a union of connected components of ${{\rm Hilb}^{G}_{h}}(G {\times}^H X)$ and
\[ {{\rm Hilb}^{G}_{h}}(G {\times}^H X)=\coprod_{i}{{\rm Hilb}^{H}_{h_{i}}}(X) \]
\end{lem}
\begin{proof}
This is an immediate consequence of the proof of \cite[Lemma 3.2]{brion}.
\end{proof}

Now we can easily prove Theorem \ref{a} as follows.

\begin{proof}[(proof of Theorem \ref{a})]
By Proposition \ref{slice}, there is a saturated open immersion $U \hookrightarrow X$ and a strongly {\'e}tale surjective morphism
$G\times^{G_x} V \twoheadrightarrow U$.
Applying Lemma \ref{satop} and \ref{stet} to these morphisms, we get
\[
{{\rm Hilb}^{G}_{h}}(G\times^{G_x} V) \twoheadrightarrow {{\rm Hilb}^{G}_{h}}(U)
\hookrightarrow {{\rm Hilb}^{G}_{h}}(X),
\]
where the first morphism is {\'e}tale surjective and the second is an open immersion.
Now the first assertion follows from Lemma \ref{ind}.
By applying Lemma \ref{stet} to the strongly {\'e}tale morphism $V \longrightarrow T_{x}V$ of Proposition \ref{smslice},
the second assertion follows.
\end{proof}

The proof of Theorem \ref{b} is slightly complicated because we have to chase elements
of the invariant Hilbert schemes.

\begin{proof}[(proof of Theorem \ref{b})]
We have the following commutative diagram.
\[
\xymatrix{
&Z \ar[ld]_{j} \ar[d]^i\\
G\times^{G_x}V \ar[r] \ar[d] & X \ar[d]\\
(G\times^{G_x}V)\cquot G \ar[r] & X\cquot G
}
\]
Here $i$ denotes the inclusion.
The morphism $j$ is defined as follows.
By assumption, the scheme theoretic image of $Z \longrightarrow X\cquot G$ is a reduced closed point.
Let $y \in (G\times^{G_x}V)\cquot G$ be the image of $(e,x) \in G \times V$.
Then $j$ is the morphism induced by the constant map $Z \longrightarrow y\in  (G\times^{G_x}V)\cquot G$ and $i \colon Z \longrightarrow X$.
Note that the square in the above diagram is cartesian.
Since $i$ is a closed immersion, so is $j$.

Next we consider the following diagram
\begin{align}\label{key diagram}
\xymatrix{
Z \ar[r]_j \ar@/^5mm/[rr]^{f} & G\times^{G_x}V \ar[r]_p & G/G_x\\
W \ar[u] \ar[r] & V \ar[u] \ar[r] & \overline{e}. \ar[u]
}
\end{align}

Here $p$ is the morphism induced by the first projection $G\times V \longrightarrow G$
and the two squares are cartesian. Note that $f$ gives a homeomorphism
of underlying spaces and in particular is a quasi-finite morphism.
Therefore $W$ is finite over $k$ and its Hilbert function $h'$ has finite support.
By the construction we have
\[
\xymatrix@C=36pt@R=5pt{
{{\rm Hilb}^{G_x}_{h'}}(V) \ar[r] & {{\rm Hilb}^{G}_{h}}(G\times^{G_x}V) \ar[r] & {{\rm Hilb}^{G}_{h}}(X) \\
(W \longrightarrow V) \ar@{}[u]|{\rotatebox[origin=c]{90}{$\in$}} \ar@{|->}[r] &
{(j\colon Z\longrightarrow G\times^{G_x}V)} \ar@{}[u]|{\rotatebox[origin=c]{90}{$\in$}} \ar@{|->}[r] &
{(i\colon Z\longrightarrow X)} \ar@{}[u]|{\rotatebox[origin=c]{90}{$\in$}}
}
\]
where the first morphism is the inclusion in Lemma \ref{ind}
and the second one is {\'e}tale by Lemma \ref{stet}.
The first part of the theorem follows.

The proof of the second assertion is similar to that of Theorem \ref{a}.
\end{proof}

\begin{rmk}\label{red}
We use the notation of the proof. If we put $n=\sum_{M\in {\rm Irr}(G_x)} h'(M)$ (which is well defined because $h'$ is
support finite), we have the following commutative diagram.
\[
\xymatrix{
{\rm Hilb}^{G_x}_{h'}(V) \ar[rr] \ar[rd] & & {{\rm Hilb}_{n}(V)}\\
&{{\rm Hilb}_{n}(V)}^{G_x} \ar[ru] &
}
\]
If $X$ is smooth at $x$, then we also have the following commutative diagram.
\[
\xymatrix{
&{\rm Hilb}^{G}_{h}(X)&\\
{\rm Hilb}^{G_x}_{h'}(V) \ar[ru]^{\text{{\'e}tale}} \ar[r] \ar@/^5.5mm/[rr]^{\alpha}
\ar@/_4mm/[rrd]^{\gamma}_{\text{{\'e}tale}}&
{\rm Hilb}^{G_x}_{h'}(T_xV) \ar[r] \ar[rd]^{\beta}&
{{\rm Hilb}_{n}(T_xV)}\\
&&{{\rm Hilb}_{n}(T_xV)}^{G_x} \ar[u]
}
\]
Here 
\begin{itemize}
\item the scheme ${{\rm Hilb}_{n}(T_xV)}^{G_x}$ is the fixed points scheme;
\item the morphism $\alpha$ maps $W$ to length $n$ closed subscheme of $T_xV$ supported at $0$;
\item the morphism $\beta$ is the canonical open immersion.
\end{itemize}
Therefore $\gamma$ is an {\'e}tale morphism and
$\gamma (W)$ represents $G_x$-stable closed subscheme of $T_xV$ supported at $0$.

\end{rmk}

To prove Corollary \ref{c}, we need some facts.
The following lemma is well-known (see for example \cite[Corollary to Theorem 5.4]{fo73}).

\begin{lem}\label{fix}
Let $H$ be a reductive group over $k$, 
and $Y$ be a separated $H$-scheme of finite type over $k$.
Let $Y^H$ be the fixed points scheme.
If $Y$ is smooth at a closed point $y \in Y^H$, then $Y^H$ is
also smooth at $y$.
\end{lem}

As for the smoothness of punctual Hilbert schemes, we have the following facts.

\begin{prop}\label{punc}
Let ${\rm Hilb}_{n}({\mathbb{A}}^d)$ be the punctual Hilbert scheme of $n$ points.
\begin{enumerate}
\item[{\rm (1)}] If $d \leq 2$, then ${\rm Hilb}_{n}({\mathbb{A}}^d)$ is smooth {\rm \cite[Theorem 2.4]{fo68}}.
\item[{\rm (2)}]  If $d \leq 3$ and $Z\in {\rm Hilb}_{n}({\mathbb{A}}^d)(k)$
is Gorenstein as a scheme, then ${\rm Hilb}_{n}({\mathbb{A}}^d)$ is smooth at $Z$ {\rm \cite[Corollary 2.6]{cn}}.
\end{enumerate}
\end{prop}

\begin{proof}[(proof of Corollary \ref{c})]
Take a smooth slice $V$ at $x$. By Proposition \ref{smslice},
${\rm dim}V={\rm dim}_xX-{\rm dim}Gx$.
We use the notation from the proof of Theorem \ref{b} and Remark \ref{red}.
By Theorem \ref{b} and its proof, ${{\rm Hilb}^{G}_{h}}(X)$ is smooth at $Z$ if and only if
${\rm Hilb}^{G_x}_{h'}(V)$ is smooth at $W$. 

(1) By Proposition \ref{punc} (1), ${\rm Hilb}_{n}(T_xV)$ is smooth.
By Lemma \ref{fix}, ${{\rm Hilb}_{n}(T_xV)}^{G_x}$ is smooth.
Since $\gamma$ is {\'e}tale, ${\rm Hilb}^{G_x}_{h'}(V)$ is also smooth.

(2)
Consider the following diagram, where the square is in the proof of Theorem \ref{b} (diagram (\ref{key diagram})).
\[
\xymatrix{
G/G_x \ar[r] \ar@/^4mm/[rr]^{{\rm id}}&Z \ar[r]_f &G/G_x\\
&W=f^{-1}(\overline{e}) \ar[u] \ar[r]&\overline{e} \ar[u]
}
\]
Here the morphism $G/G_x \longrightarrow Z$ is the natural morphism from the reduced scheme of $Z$.
Therefore every local ring of $W$ is the quotient of a local ring of $Z$
by a regular sequence. Hence $W$ is Gorenstein if (and only if) $Z$ is Gorenstein.
By Proposition \ref{punc} (2), ${\rm Hilb}^{G_x}_{h'}(T_xV)$ is smooth at $W$
and ${\rm Hilb}^{G_x}_{h'}(T_xV)^{G_x}$ is also smooth at $W$.
Now the statement follows from Remark \ref{red}.

(3) If $Z$ is reduced, then $W=\Spec k$ and $h'$ is trivial.
Therefore ${\rm Hilb}^{G_x}_{h'}(V)=V^{G_x}$ and this is smooth by Lemma \ref{fix}. 

\end{proof}

\section{An example}\label{an example}

Set $G= {\mathbb{G}}_{{\rm m}}^{2}$ and $X= {\mathbb{A}}^{4}$.
Define $G \curvearrowright X$ by $(s,t)\cdot (x,y,z,w)=(s^{2}x,s^{4}y,s^{-8}z,s^{-3}tw)$.
Let $h \colon {\rm Irr}(G) \longrightarrow {\mathbb{Z}}$ be the Hilbert function of the orbit of
any point $(x,y,z,w)$ with $w=0$ and $xyz \neq 0$.
That is, $h$ is the Hilbert function of $G/(\{1,-1\}\times {\mathbb{G}}_{{\rm m}})$ and therefore

\begin{align*}
h \colon {\rm Irr}(G)= {\mathbb{Z}}^{2} \longrightarrow {\mathbb{Z}}; (l,l') \mapsto  
\begin{cases}
1 & \text{$l$ is even and $l'=0$}\\
0 & \text{otherwise.}
\end{cases}
\end{align*}

We claim that $Y:= {\rm Hilb}^{G}_{h}(X)$ is smooth.
Let $[Z] \in Y$ where $Z\subset X$ is a closed $G$-stable subscheme with Hilbert function $h$.
Then, there are the following three possibilities:
\begin{enumerate}
\item $Z=G\cdot P$ where $P=(x,y,z,w)$ is a point with $w=0, xyz \neq 0$;
\item $Z=G\cdot P$ where $P=(x,y,z,w)$ is a point with $y=w=0, xz \neq 0$;
\item $\Supp Z=G\cdot P$ where $P=(x,y,z,w)$ is a point with $x=w=0, yz \neq 0$ and the coordinate ring of $Z$ is the form of
$k[x,y,z]/(x^{2}, y^{2}z-a)$ with $a\in k^{*}$.
\end{enumerate}

In the cases of (1) and (2), $Y$ is smooth at $[Z]$ by Corollary \ref{c} (3).
In the case of (3), $G_{P}=\{\pm1,\pm \sqrt{-1}\}\times {\mathbb{G}}_{{\rm m}}$ and the Hilbert function $h'$ in Theorem \ref{b} is 
\begin{align*}
h' \colon {\rm Irr}(G_{P})=\{0,1,2,3\}\times {\mathbb{Z}} \longrightarrow {\mathbb{Z}}; (l,l') \mapsto 
\begin{cases}
1 & \text{$l=0,2$ and $l'=0$}\\
0 & \text{otherwise.}
\end{cases}
\end{align*}
In particular, we have $\sum_{M\in {\rm Irr}(G_{P})}h'(M)=2$.
Note that the {\'e}tale slice at $P$ is three dimensional and $ {\rm Hilb}_{2}( {\mathbb{A}}^{3})$ is smooth.
Thus, by Remark \ref{red}, $Y$ is smooth at $[Z]$.

We can explicitly construct $Y$ and the universal family.
Ideals of subschemes of type (1) are the form of $( y^{2}z-a,x^{2}-by,w,x^{4}z-ab^{2})$,
type (2) $(x^{4}z-a,y,w)$ and type (3) $(x^{2}, y^{2}z-b,w)$ with $a,b\in k^{*}$.
Consider the following flat families of $G$-subschemes of $ {\mathbb{A}}^{4}$:
\begin{align*}
&\mathcal{Z}=\Spec k[x,y,z,w,a,a^{-1},b]/(y^{2}z-a,x^{2}-by,w,x^{4}z-ab^{2})\subset {\mathbb{A}}^{4}\times \Spec k[a,a^{-1},b];\\
&\mathcal{W}=\Spec k[x,y,z,w,c,d,d^{-1}]/(y^{2}z-c^{2}d,cx^{2}-y,w,x^{4}z-d)\subset {\mathbb{A}}^{4}\times \Spec k[c,d,d^{-1}].
\end{align*}
Then, $ \mathcal{Z}$ is a family that consists of all type (1) and type (3) subschemes and
$ \mathcal{W}$ that of all type (1) and type (2) subschemes.
Since $Y$ is smooth, it is isomorphic to the scheme obtained by glueing $ {\mathbb{G}}_{{\rm m}}\times {\mathbb{A}}^{1}$ and
$ {\mathbb{A}}^{1}\times {\mathbb{G}}_{{\rm m}}$ by the following rule:
\[
\xymatrix@C=36pt@R=7pt{
{\mathbb{G}}_{{\rm m}}\times {\mathbb{A}}^{1} \ar@{}[r] | {\supset}& {\mathbb{G}}_{{\rm m}}^{2} \ar@{}[r]|{\simeq} &
{\mathbb{G}}_{{\rm m}}^{2} \ar@{}[r]|{\subset}&{\mathbb{A}}^{1} \times {\mathbb{G}}_{{\rm m}}\\
&(a,b)\ar@{}[u]|{\rotatebox[origin=c]{90}{$\in$}} \ar@{|->}[r] &(1/b,ab^{2}) \ar@{}[u]|{\rotatebox[origin=c]{90}{$\in$}}.&
}
\]

\begin{ack}
The author would like to thank Professor Tomohide Terasoma for his advice and comments.
The author would like to thank the referee for useful suggestions that helped
him to improve the original manuscript.
The author is supported by the Program for Leading Graduate Schools, MEXT, Japan.
\end{ack}

\end{document}